\documentclass[a4paper,oneside,10pt]{amsart}

\usepackage{amsmath,amssymb,amsthm,latexsym,array,amsfonts,eufrak,mathrsfs,nomencl,hyperref,threeparttable,textcomp}
\usepackage{graphicx,setspace,subfigure,epigraph,color,booktabs,caption,multirow,titlesec}
\usepackage[utf8]{inputenc}
\usepackage[english]{babel}
\usepackage[a4paper,top=3cm,bottom=3cm,left=3.1cm,right=3.1cm]{geometry}
\numberwithin{equation}{section}
\titleformat{\section}[hang]{\normalfont\fillast}{\thesection. }{0pt}{\scshape}
\titleformat{\subsection}[hang]{\normalfont\fillast}{\thesubsection. }{0pt}{\scshape}
\titleformat{\subsubsection}[hang]{\normalfont\fillast}{\thesubsubsection. }{0pt}{\scshape}

\title{Nonnegative polynomials and their Carathéodory number}
\author{Simone Naldi}

\address{Università di Firenze, Dipartimento di Matematica ``Ulisse Dini''}
\address{Viale Morgagni 67/A, 50134, Firenze, Italy}
\address{mail: naldi@laas.fr}
\date{}

\newtheorem{teo}{Theorem}[section]
\newtheorem{prop}[teo]{Proposition}
\newtheorem{coroll}[teo]{Corollary}
\newtheorem{lemma}[teo]{Lemma}

\newcommand{\Ree}{\text{Re}}
\newcommand{\Imm}{\text{Im}}

\theoremstyle{definition}
\newtheorem{defin}{Definition}[section]
\newtheorem{example}{Example}[section]
\newtheorem{oss}{Remark}[section]

\begin{document}

\maketitle

\begin{abstract}
\noindent In 1888 Hilbert showed that every nonnegative homogeneous
polynomial with real coefficients of degree $2d$ in $n$ variables
is a sum of squares if and only if $d=1$ (quadratic forms), $n=2$
(binary forms) or $(n,d)=(3,2)$ (ternary quartics). In these cases,
it is interesting to compute canonical expressions for these decompositions.
Starting from Carathéodory's Theorem, we compute the Carathéodory number
of Hilbert cones of quadratic forms and binary forms.
\end{abstract}

\thispagestyle{empty}

\section*{Introduction}
The decomposability of a given nonnegative polynomial as a sum of squares
is an old and interesting problem in mathematics, with countless applications
in many research fields. It turns out that not all real polynomials, which
are nonnegative over $\mathbb{R}^n$, can also be written as the sum of squares
of other polynomials: an expression like this could be a direct certificate of
nonnegativity for the given polynomial. Hilbert, in 1888 \cite{hilbert1888darstellung},
solved this problem by means of a celebrated result, showing that every
nonnegative homogeneous polynomial $F\in{\mathbb{R}[x_1,...,x_n]}$ of degree $2d$
is a sum of squares if and only if $n=2$ (binary forms), $d=1$ (quadratic forms)
or $(n,d)=(3,2)$ (quartics in three variables).
In the 20th century many examples of nonnegative polynomials that can't be
written as a sum of squares were produced \cite{motzkin1967arithmetic,
reznick1992sums,blekherman2006there}. \smallskip \smallskip

We start from Hilbert's classical result, studying some aspect of particular decompositions.
Remark that the set of nonnegative polynomials is a convex cone in a finite-dimensional vector space.
So it is interesting to give a description of its \textit{extreme points} (or \textit{extreme rays}),
that are points that cannot be non-trivially decomposed as sum of other elements of the set.
Extreme points form a fundamental subset of the cone: by Carathéodory's Theorem (see Theorem
\ref{carat}), if the dimension of the vector space is finite and under some suitable conditions,
every element of the cone can be written as a finite sum of extreme points. Furthermore, fixed $x$ in
the cone, if we call the \emph{length} of $x$ the minimum integer $t$ such that $x$ is a sum of
$t$ extreme points, from Carathéodory's Theorem one obtains that the length of $x$ is bounded
by the dimension of the cone. Since every element of the cone $\mathbf{P}_{n,2d}$ of nonnegative
polynomials is a finite sum of extreme points, it is interesting to compute the maximum length
attained in the set $\mathbf{P}_{n,2d}$: this value is the \emph{Carathéodory number} of
$\mathbf{P}_{n,2d}$. \smallskip \smallskip

\textbf{Motivations}. Computing the Carathéodory number $\mathcal{C}(X)$ for a given convex set $X$ is a general problem
which finds applications in many fields of mathematics. Although Carathéodory's Theorem is a
classical result in convex analysis, nowadays there exist various open problems concerning $\mathcal{C}(X)$,
both for compact and for non-compact cases. This problem is also naturally related with
problems in operational research; for compact convex sets like the \textit{orbitopes} \cite{sanyal2009orbitopes}
this number is known only in some particular cases. Independently from its formal definition, it is
possible to exploit various results about Carathéodory number in numerical methods for the resolution
of optimization problems (for an example in chemical applications we refer to \cite{gardner2005reconstruction}).
For example, if you want to represent, with the method of least squares, an element $x$ of a given cone
$C$ such that $x=\sum_{i=1}^{t}{e_i}$ with respect to the class $\{e_i\}$ of extreme points of $C$,
it is clear that the degrees of freedom of this representation increase as $t$ increases. So it is interesting
to find the optimal (i.e. the minimum) $t$ such that this representation is possible for every $x$ in the cone. 
Finally, for linear programming problems, it is known that minima and maxima lie in
in the set of extreme points of the feasible polyhedron: by this, it is of first importance to characterize this set. \smallskip \smallskip 

\textbf{Main results}. The main result of this paper concerns decompositions of nonnegative binary forms. It states that the Carathéodory
number of $\mathbf{P}_{2,2d}$ is $2$ for every $d$. This means that every nonnegative binary form is generated
by a pair of extreme points, which is the statement of Corollary \ref{main}.
We mainly use the notation and some results of \cite{reznick1992sums} for the cone of nonnegative polynomials
and for the duality relations of Hilbert cones. For general results about the structure of these cones we refer to
\cite{blekherman2012nonnegative,blekherman2012algebraic,blekherman2004convexity,blekherman2006there,choi1977extremal,reznick2000some,reznick2010length,powers1999notes} and to their references. Finally, we would like to recall the work \cite{karlin1952geometry} by Karlin and Shapley, where they obtained
a similar result studying the connection between the structure of the cone of univariate polynomials that are nonnegative over compact
intervals of the real line and the moments problem (after generalized in \cite{karlin1963representation} by Karlin to the whole line). Our result
presents a different and constructive viewpoint of these ideas in the context of sums-of-squares representations of nonnegative polynomials. \\

\textbf{Acknowledgements.} The author thanks Giorgio Ottaviani and Marco Longinetti for their fundamental aid and their
constant encouragement to improve this work. The author also thanks the anonymous reviewers for the time they spent reading
this work and for all their precious advice.

\section{Preliminaries}
For the algebraic-geometric dictionary we use in this paper we refer to \cite{harris1992algebraic} and \cite{hartshorne1977algebraic},
and for a general introduction to convexity to \cite{fenchel1953convex}.
Let $\mathbb{R}$ be the field of real numbers, and let $H$ be an $\mathbb{R}$-vector space of finite dimension.
\begin{defin}
A subset $C\subseteq{H}$ is a \textit{cone} if $\text{for each} \ a\in{C} \ \text{and each} \ \alpha>0 \ \text{then} \ \alpha{a}\in{C}.$
\end{defin}
The \textit{dimension} of a cone is the dimension of its affine hull.
We are interested in \textit{closed} and \textit{convex} cones of finite dimension, and full-dimensional in the
vector space where they are defined. It is easy to see that a cone $C$ is convex if and only if, given $a,b\in{C}$, $a+b\in{C}$.
Now, if $C$ is a closed cone, an element $e\in{C}$ is an \textit{extreme point} for $C$ if every decomposition $e=f_1+f_2$,
with $f_i\in{C}$, implies $f_i=c_ie$, for some $c_i\geq{0}$. So, a point $x$ in a cone is extreme if and only if it cannot
be non-trivially decomposed as the sum of two other elements of the cone. 

We denote by $\text{Ext}(C)$ the set of extreme points of $C$. It is clear that $\mu\cdot\text{Ext}(C)=\text{Ext}(C)$
for every $\mu\in\mathbb{R}^+$. Moreover, in a closed convex cone the set of extreme points is a subset of the boundary of the cone.
We refer to the following two classical facts: the first is a version for convex cones over a compact set
of Krein-Milman Theorem \cite{fenchel1953convex}, for the second the reader can find an easy proof on \cite{reznick1992sums}.
\begin{teo}[Krein-Milman]\label{milm}
Let $C$ be a convex cone such that the following holds: if $0\neq x\in C$, then $-x\notin C$. Then $C$ is the convex hull
of its extreme points.
\end{teo}
\begin{teo}[Carathéodory]\label{carat}
Let $C$ be a convex cone of dimension $N$ containing $0$, and $z=\sum_{j=1}^{b}{x_j}$ with $x_j\in{C}$.
Then there exist $\{y_1,\dots,y_N\}\subseteq{C}$ such that $z=\sum_{i=1}^{N}{y_i}$ and that $\forall i \ \exists j$ s.t.
$y_i=\epsilon_j{x_j}$, for some $\epsilon_j\leq{1}$.
\end{teo}
So one obtains
\begin{coroll}\label{corollarioestremali}
Let $C$ be a closed convex cone such that the hypothesis of Theorem \ref{milm} are satisfied. Let $N$ be the dimension of $C$.
Then every $x\in{C}$ has a representation as a sum of at most $N$ extreme points of $C$.
\begin{proof}
Let $x\in{C}$. Then, by Theorem \ref{milm}, $x$ can be written as sum of a finite number of extreme points of $C$, and by Theorem \ref{carat}
this number can be reduced to $N$.
\end{proof}
\end{coroll}
We call \textit{extremal} a decomposition of an element $x\in C$ of type $x=\sum e_i$ for $e_i\in\text{Ext}(C)$. By Corollary \ref{corollarioestremali}
we can give the following definition.
\begin{defin}
Let $C$ be a closed convex cone of finite dimension such that if $0\neq x\in C$, then $-x\notin C$, and let $x \in C$. We call the \textit{length} of $x$
the integer value $$\text{h}(x)=\min \left\{r\in\mathbb{N} \ \Big| \ x = \sum_{i=1}^r{e_i} \ \text{for some} \ e_i\in\text{Ext}(C)\right\}$$ and the
\textit{Carathéodory number} of the cone $C$ the value $$\mathcal{C}(C)=\max_{x\in C}{\Big\{\text{h}(x)\Big\}}.$$
\end{defin}

\begin{oss}\label{CC}
$\text{h}(x)\leq \dim(C)$ for every $x\in C$, and so $\mathcal{C}(C)\leq \dim(C)$.
\end{oss}

\section{The Carathéodory number for Hilbert cones}
We will denote $\mathbf{F}_{n,2d}$ the vector space $\mathbb{R}[x_1,...,x_n]_{2d}$ of real homogeneous
polynomials of degree $2d$ in $n$ variables, which is the $2d-$th homogeneous component of the gradued
ring $\mathbb{R}[x_1,...,x_n]$. Then, by
\begin{equation}
\mathbf{P}_{n,2d}=\Big\{P\in\mathbf{F}_{n,2d} \ \Big| \ P(X)\geq 0, \ \forall \ X=(x_1,\dots,x_n)\in\mathbb{R}^n\Big\} 
\end{equation}
the set of nonnegative polynomials in $n$ variables of degree $d$, and by
\begin{equation}
\mathbf{\Sigma}_{n,2d}=\Big\{\sum_{i=1}^t{g_i^2} \ \Big| \ g_i\in\mathbf{F}_{n,d}, t\in\mathbb{N}\Big\}
\end{equation}
the set of sums of squares. Of course one has the trivial inclusion $\mathbf{\Sigma}_{n,2d}\subseteq\mathbf{P}_{n,2d}$ for
each $(n,d)$. In 1888, Hilbert characterized the cases when this inclusion is an equality.
\begin{teo}[Hilbert, 1888]\label{hilbert}
$\mathbf{P}_{n,2d}=\mathbf{\Sigma}_{n,2d}$ if and only if one of the following cases occurs:
\begin{itemize}
\item[(a)] \ \ $n=2$ \ $\mathrm{(binary \ forms)}$
\item[(b)] \ \ $d=1$ \ $\mathrm{(quadratic \ forms)}$
\item[(c)] \ \ $(n,d)=(3,2)$ \ $\mathrm{(ternary \ quartics)}$
\end{itemize}
\end{teo}
If we denote by 
\begin{equation}\label{psi}
\Psi=\Big\{(n,1),(2,d),(3,2) \ \Big| \ n,d\in\mathbb{N}\Big\},
\end{equation}
Hilbert's theorem says that $\mathbf{P}_{n,2d}=\mathbf{\Sigma}_{n,2d}$ if and only if $(n,d)\in\Psi$.
It is easy to prove the following proposition (for a proof, see \cite[Prop.3.6]{reznick1992sums}).
\begin{prop}
 $\mathbf{P}_{n,2d}$ and $\mathbf{\Sigma}_{n,2d}$ are full-dimensional closed convex cones of $\mathbf{F}_{n,2d}$.
\end{prop}
$\mathbf{P}_{n,2d}$ and $\mathbf{\Sigma}_{n,2d}$ are the so-called \textit{Hilbert cones} of polynomials.
For $F\in\mathbf{F}_{n,2d}$ denote $V(F)$ the complex algebraic variety encoding the points where
the polynomial $F$ vanishes, and $\text{V}_{\mathbb{R}}(F)$ the set of real points of $V(F)$.
\begin{prop}
If $F\in\partial\mathbf{P}_{n,2d}$ then every point of $\mathrm{V}_{\mathbb{R}}(F)\subset\mathbb{R}^n$ is singular.
\begin{proof}
It is clear
that the interior of $\mathbf{P}_{n,2d}$ consists of all polynomials that are positive definite over $\mathbb{R}^n$, that is, for
every $F$ in the interior, $F(X)>0$ for every $X\in\mathbb{R}^n$, $X\neq 0$. Let $F\in\partial\mathbf{P}_{n,2d}$. If $a\in\mathbb{R}^n$
such that $F(a)=0$, then by nonnegativity $\nabla F(a)=0$, that is $a$ is a singular point for $\mathrm{V}_{\mathbb{R}}(F)$.
\end{proof}
\end{prop}
\begin{prop}\label{estre}
If $(n,d)\in\Psi$ then $\mathrm{Ext}(\mathbf{P}_{n,2d})\subseteq\mathbf{F}^2_{n,d}=\Big\{P^2 \ | \ P\in\mathbf{F}_{n,d}\Big\}$.
\begin{proof}
In fact if $(n,d)\in\Psi$, then $\mathbf{P}_{n,2d}=\mathbf{\Sigma}_{n,2d}$. So, let $F\in\text{Ext}(\mathbf{P}_{n,2d})$,
then $F=\sum_{k=1}^{w}{g_k^2}$ where we can choose $w$ as the smallest integer with this property. So $g_i\neq \alpha g_j$
for $i\neq j$ for every $\alpha$, because otherwise we could write
\begin{equation}
 F=\sum_{k\neq i,j}{g_k^2}+(1+\alpha^2)g_j^2,
\end{equation}
which is a sum of $w-1$ squares, and this is a contradiction. Since $F$ is extreme, the only possibility is that $w=1$, so that $F$ is a square.
\end{proof}
\end{prop}
So, if $(n,d)\in\Psi$, then $\mathrm{Ext}(\mathbf{P}_{n,2d})\subseteq\Big[\mathbf{F}^2_{n,d}\cap\partial\mathbf{P}_{n,2d}\Big]$.
In this paper we address the problem of characterizing the subset of extreme points of $\mathbf{P}_{n,2d}$
and in calculating in some cases the Carathéodory number of this cone. We remember that, by Remark \ref{CC} one obtains the upper bound
\begin{equation}
\mathcal{C}(\mathbf{P}_{n,2d})\leq\binom{n+2d-1}{n-1}=\dim\mathbf{P}_{n,2d}=\dim\mathbf{F}_{n,2d}. 
\end{equation}

\subsection{Apolarity and duality}
Let us consider for every $d$ the set of polynomials
$$\mathbf{Q}_{n,2d}=\Big\{F=\sum_{k=1}^{r}{(\alpha_k\cdot X)^{2d}} \ \Big| \ \alpha_k{\in{\mathbb{R}^n}}, r\in\mathbb{N}\Big\}$$
(where $\alpha_k\cdot X=\alpha_1x_1+\dots+\alpha_nx_n$) of finite sums of $2d-$th powers of linear forms, which is a closed
convex subcone of $\mathbf{\Sigma}_{n,2d}$. Let $\mathbb{P}^{\ell}$ denote the $\ell-$dimensional projective space over the field $\mathbb{R}$.
Since every element of $\mathbf{Q}_{n,2d}$ is a sum of $2d-$th powers, it is clear that the subset of extreme points of $\mathbf{Q}_{n,2d}$
consists of all $2d-$th powers of linear real forms: this set is strictly linked to the image $\mathscr{V}_{n-1,2d}(\mathbb{P}^{n-1})$
of the Veronese map, defined by
\begin{equation}
 \mathscr{V}_{n-1,2d} \colon \mathbb{P}^{n-1} \to \mathbb{P}^{\binom{n+2d-1}{n-1}-1}\cong{\mathbb{P}(\text{S}^{2d}(\mathbb{R}^{n}))},
\end{equation}
with $\mathscr{V}_{a,b}([x_0:\dots:x_a])=[\cdots:X^I:\cdots]$,
where $I$ ranges over the set $\mathcal{I}(a,b)=\{I=(i_0,\dots,i_a) \ | \ i_j\in\mathbb{N} \ \mathrm{and} \ \sum{i_j}=b\}$.
This image parametrizes the variety of $2d-$th powers of linear forms, subvariety of the space $\mathbb{P}^{\binom{n+2d-1}{n-1}-1}$
of hypersurfaces of $\mathbb{P}^{n-1}$ of degree $2d$: so the Zariski closure $\overline{\text{Ext}(\mathbf{Q}_{n,2d})}$
coincides with the Veronese variety in the projective space $\mathbb{P}^{\binom{n+2d-1}{n-1}-1}$. We obtain that the cone $\mathbf{Q}_{n,2d}$
is the convex hull of the Veronese variety, that is the so-called \textit{Veronese orbitope} \cite{sanyal2009orbitopes}. \smallskip \smallskip 

Now, let $\mathbb{R}[\partial_{x_1},\dots,\partial_{x_n}]_{2d}=\mathbf{F}_{n,2d}'$ be the dual ring of $\mathbf{F}_{n,2d}$. We have that
$\mathbf{F}_{n,2d}$ is generated by the set of monomials $\{X^I\}_{I\in\mathcal{I}(n,2d)}$ where $X^I=x_1^{i_1}\cdots x_n^{i_n}$
while $\mathbf{F}_{n,2d}'$ is generated by the set of monomials of formal derivatives $\{\frac{\partial_I}{c(I)}\}_{I\in\mathcal{I}(n,2d)}$ where
$\partial_I=\frac{1}{\prod{i_j!}}\partial_{x_1}^{i_1}\cdots \partial_{x_n}^{i_n}$ and $c(I)=\frac{(2d)!}{\prod{i_j!}}$.
Consider the natural map $\sigma \colon \mathbf{F}_{n,2d} \to \mathbf{F}_{n,2d}'$ induced by the choice of these bases. One has
$\sigma(X^I)=\frac{\partial_I}{c(I)}$ for all $I$. This map defines the following bilinear form on $\mathbf{F}_{n,2d}$:
the map $\tau \colon \mathbf{F}_{n,2d}\times\mathbf{F}_{n,2d} \to \mathbb{R}$ such that $\tau(F,G)=\left(\sigma(G)\right)(F)$.
Using coordinates, if $F=\sum c(I)a_F(I)X^I$ and $G=\sum c(I)a_G(I)X^I$ we obtain
$$\tau(F,G)=\sum c(I)a_F(I)a_G(I)=\tau(G,F).$$
Under this product, $\mathbf{P}_{n,2d}$ and $\mathbf{Q}_{n,2d}$ are mutually dual \cite{reznick1992sums}.
Now, fixed $F\in\mathbf{F}_{n,2d}$ and $1\leq i\leq 2d-1$, consider the $i$-th \textit{apolarity map}
\begin{equation}
\mathrm{AP}_F(i,2d-i;n) \colon \mathbf{F}_{n,2d-i} \to \mathbf{F}_{n,i}
\end{equation}
sending $G \mapsto \left(\sigma(G)\right)(F)$. $\mathrm{AP}_F(i,2d-i;n)$ is linear: the matrix of the map $\mathrm{AP}_F(d,d;n)$
is the \textit{catalecticant matrix} of $F$, and we denote it by $\mathrm{H}_F$. If $F$ is a quadratic form of matrix $M$,
then $\mathrm{H}_F=M$; if $F$ is a binary form of degree $2d$, $\mathrm{H}_F$ is a Hankel matrix of order $d+1$ whose terms
on the diagonals are the coefficients of $F$. The catalecticant of $F\in\mathbf{P}_{n,2d}$ is strictly linked to the decomposability
of $F$ as a sum of squares or powers of linear forms, a special instance of the Waring Problem for polynomials. For example,
$\mathbf{\Sigma}_{n,2d}^*$ (the convex cone dual to $\mathbf{\Sigma}_{n,2d}$) is the set of nonnegative polynomials whose catalecticant
matrix is positive semidefinite; so, if $(n,d)\in\Psi$, $\mathbf{Q}_{n,2d}=\mathbf{\Sigma}_{n,2d}^*$ has this property. 
In \cite[Theorem 4.6]{reznick1992sums}, the author shows the following interesting fact:
\begin{teo}[Reznick]
For every $F\in\mathbf{Q}_{n,2d}$ let $w(F)$ be its length in $\mathbf{Q}_{n,2d}$. Then, for every $F$,
$w(F)\geq{\mathrm{rk}(\mathrm{H}_F)}$. The equality holds for every $F\in\mathbf{Q}_{n,2d}$ if and only if $(n,d)\in\Psi$.
\end{teo}
This result, considered in the context of Carathéodory number means that $\mathcal{C}(\mathbf{Q}_{n,2d})$
is always greater than the maximum rank of the catalecticant matrix of a generic sum of $2d-$th powers of linear forms on $n$ variables.
\begin{prop} We obtain the following results about $\mathcal{C}(\mathbf{Q}_{n,2d})$:
\begin{table}[h]
\begin{center}
\begin{tabular}{cc}
\toprule[1.2pt]
$(n,d)$ & $\mathcal{C}(\mathbf{Q}_{n,2d})$\\
\midrule[1pt]
$(n,1)$ & $\mathcal{C}=n$\\
%\midrule
$(2,d)$ & $\mathcal{C}=d+1$\\
%\midrule
$(3,2)$ & $\mathcal{C}=6$\\
\midrule[1pt]
 $\notin\Psi$ & $\binom{n+d-1}{n-1}\leq\mathcal{C}\leq\binom{n+2d-1}{n-1}$\\
\bottomrule[1.2pt]
\end{tabular}
\end{center}
\end{table}
\end{prop}

\subsection{Quadratic forms}
Let us consider the cones of nonnegative quadratic forms in $n$ variables, that is
$\mathbf{P}_{n,2}=\{F\in\mathbf{F}_{n,2} \ | \ F(X)\geq 0 \ \forall \ X\in\mathbb{R}^n\}$, $n\in\mathbb{N}$.
There is a natural identification between $\mathbf{P}_{n,2}$ and the cone of $n\times{n}$ symmetric positive semidefinite matrices
with real entries: every nonnegative quadratic form $F\in\mathbf{P}_{n,2}$ has a representation
$F(X)={}^tM_FX M_F, \ \text{for every} \ X\in\mathbb{R}^n,$ with $M_F$ positive semidefinite.
Diagonalizing $M_F$ one obtains a canonical form for $F$, that is
\begin{equation} \label{quad}
 F(X)=\sum_{i=1}^{\text{rk} M_F}{(t_i\cdot X)^2},
\end{equation}
where $\text{rk} M_F$ is the rank of the matrix.
\begin{prop}
 $\mathrm{Ext}(\mathbf{P}_{n,2})=\mathbf{F}_{n,1}^2=\partial\mathbf{P}_{n,2}, \ \forall \ n\in\mathbb{N}$.
\begin{proof}
We already know that $\mathrm{Ext}(\mathbf{P}_{n,2})\subseteq\mathbf{F}_{n,1}^2$ by Proposition \ref{estre} and
$\mathbf{F}_{n,1}^2=\partial\mathbf{P}_{n,2}$ because every polynomial of the boundary has at least one double real root
and the bound on the degree necessarily forces the equality. Let now $F\in\mathbf{F}_{n,1}^2$.
Then $F=(\alpha\cdot X)^2=(\alpha_1 x_1+\dots+\alpha_n x_n)^2$ where
 $\alpha=(\alpha_1\dots\alpha_n)\in\mathbb{R}^n$. The (real) zeros of $F$ compose the set of points of $\mathbb{R}^n$
orthogonal to $\alpha$, that is the hyperplane $\alpha^{\bot}$.
If $F=F_1+F_2$ with $F_i\in\mathbf{P}_{n,2}$, we obtain $\alpha^{\bot}=V(F_1)\cap V(F_2)$, that is necessarily $V(F_1)=V(F_2)=\alpha^{\bot}$.
So there exist $\lambda_1,\lambda_2$ s.t. $F_i=\lambda_i F$ and so $F\in\mathrm{Ext}(\mathbf{P}_{n,2})$.
\end{proof}
\end{prop}
Now, since every square of a linear form is an extreme point of the cone $\mathbf{P}_{n,2}$, the decomposition \eqref{quad} is extremal
and it is the minimal decomposition of $F$ as a sum of squares of linear forms. That is, $\text{h}(F)=\text{rk}M_F$
for every $F\in\mathbf{P}_{n,2}$ and so we conclude that the Carathéodory number $\mathcal{C}(\mathbf{P}_{n,2})$ is the maximum rank of
a nonnegative quadratic form on $n$ variables, that is
\begin{teo}
 $\mathcal{C}(\mathbf{P}_{n,2})=n$, for every $n\in\mathbb{N}$.
\end{teo}
We remark in the following example that there is no possibility of uniqueness of the extremal representations in $\mathbf{P}_{n,2}$,
and that generically speaking the set of extremal representations can be large. 
\begin{example} Let $G_n(X)=x_1^2+\dots+x_n^2$, $n\geq 1$, here represented with respect to the orthonormal canonical base of $\mathbb{R}^n$.
If we impose
\begin{align*}
G_n(x_1,\dots,x_n) & = x_1^2+\dots+x_n^2=\sum_{k=1}^{n} \Big(\alpha_{k1}x_1+...+\alpha_{kn}x_n\Big)^2=\\
 & = \sum_{k=1}^{n} \Big(\alpha_{k1}^2x_1^2+...+\alpha_{kn}^2x_n^2+2\sum_{i\lneq j}\alpha_{ki}\alpha_{kj}x_ix_j\Big),
\end{align*}
we obtain $n+\binom{n}{2}=\frac{n^2+n}{2}$ quadratic conditions
on the matrix $M=(\alpha_{ij})$, that force $M\in{\text{SO}(n,\mathbb{R})}$. So the family of extremal representations of $G_n$ is parametrized by
the group $\text{SO}(n,\mathbb{R})$.
\end{example}

\subsection{Binary forms}
A polynomial $F\in\mathbf{F}_{2,2d}$ can be typically expressed in the following form:
$$F(x,y)=\displaystyle\sum_{j=0}^{d}c_jx^jy^{d-j},$$
for some real numbers $c_j$. We suppose now that $F$ is positive over $\mathbb{R}^2$.
By the Fundamental Theorem of Algebra, over the complex field $F$ splits into a product of linear forms. By nonnegativity, the
real roots of $F$ have even multiplicity, and so
\begin{equation}\label{bin2}
F(x,y)=\prod_t{(a_tx-b_ty)^{2k_t}}\prod_{\alpha}{(a_{\alpha}x-b_{\alpha}y)(\overline{a}_{\alpha}x-\overline{b}_{\alpha}{y})}
\end{equation}
where the first product ranges over the real roots of $F$ (and the multiplicity of the $t$-th root is $2k_t$), while the second
one ranges over the set of couples of the complex roots with their conjugates. So $F$ is the product of a square times a product
of two complex conjugated polynomials, that is $F=R^2(C_1+iC_2)(C_1-iC_2)=(RC_1)^2+(RC_2)^2$.
Then we can see first of all that every nonnegative binary form is always a sum of at most two squares.

%%%%%%%%%%%%%%%%%%%%%%%%%%%%%%%%%%%%%%%%%%%%%%%%%%%%%%%%%%%%%%%%%%% CORRETTO FINO A QUI %%%%%%%%%%%%%%%%%%%%%%%%%%%%%%%%%%%%%%%%%%%%%%%%%%%%%%%%%%%%%%%%%%%%%%%%%%%%%

\begin{teo} Let $F\in\mathbf{P}_{2,2d}$. Then $F$ is extreme for the cone $\mathbf{P}_{2,2d}$ if and only if $F$
is the square of a polynomial with only real roots.
\begin{proof}
First suppose that $F=G^2$ has only real roots and $F=F_1+F_2$ with $F_i\in{\mathbf{P}_{2,2d}}$. One has $V(F)=V(F_1)\cap{V(F_2)}$
by nonnegativity, so $V(F)\subseteq{V(F_i)}$ for $i=1,2$. Moreover, every linear form of the factorization of $F$ must appear also in the factorization of the $F_i$'s.
In fact, let $L(x,y)=ax-by$ be a linear form that divides $F$; then $0=F(b,a)=F_1(b,a)+F_2(b,a)$ and by the nonnegativity $(b,a)$ nullifies both $F_i$'s, so $L$ divides them.
Now, if $F=L^2H$, $L$ divides $F_i$, and so does $L^2$, that is $F_i=L^2J_i$. Proceeding by induction, it is easy to see that if $F=L^{2t}H$ then $L^{2t}$ divides the $F_i$'s.
Iterating this process for each linear form of the factorization of $F$, since $F,F_1,F_2$ have degree $2d$ and $F$ has only real roots, we obtain that
$V(F)=V(F_1)=V(F_2)=L_1^{2t_1}\cup\cdots\cup{L_k^{2t_k}}$
and that $F,F_1,F_2$ define the same hypersurface. So there exist $\lambda_1,\lambda_2\in\mathbb{R}$, nonnegative, such that $F_i=\lambda_i{F}$, from which $F$ must be extreme.

In the other direction, let us suppose that $F$ is extreme in $\mathbf{P}_{2,2d}$ and assume otherwise that in the irreducible factorization of $F$
two conjugates complex factors appear. In particular, $F$ is a perfect square by Proposition \ref{estre}, and so we write:
$F=Q^2\cdot[(x-(\alpha+i\beta)y)(x-(\alpha-i\beta)y)]^2$ with $\beta\neq 0$. Then
$$ F=Q^2\cdot[(x-\alpha y)^2+(\beta y)^2]^2=Q^2\cdot(x-\alpha y)^4+Q^2\cdot(\beta y)^4+2Q^2\cdot(x-\alpha y)^2(\beta y)^2,$$
which is a non-trivial decomposition of $F$ in $\mathbf{P}_{2,2d}$. This is a contradiction because $F$ is extreme by assumption.
\end{proof}
\end{teo}

So, for binary forms there exist perfect squares that are not extreme points of $\mathbf{P}_{2,2d}$: for example, all polynomials
$(x^d+y^d)^2$ belong to $\mathbf{F}^2_{2,d}\setminus\mathrm{Ext}(\mathbf{P}_{2,2d})$ for every $d\geq 2$.

\begin{defin}
Let be $F\in{\mathbf{P}_{2,2d}}$. We call a \textit{partition of the roots} of $F$, a pair $(A,\overline{A})\in\mathbf{F}_{2,d}\times\mathbf{F}_{2,d}$
of conjugate polynomials such that
\begin{enumerate}
\item $F=A\cdot\overline{A}$
\item for all $\alpha$ s.t. $F(\alpha)=0$, then $A(\alpha)=0$ if and only if is $\overline{A}(\overline{\alpha})=0$.
\end{enumerate}
\end{defin}
So, if $F\in\mathbf{P}_{2,2d}$, and $(A,\overline{A})$ is a partition of the roots of $F$, with $A(x,y)=G(x,y)+i\cdot H(x,y)$ and
$\overline{A}(x,y)=G(x,y)-i\cdot H(x,y)$, one obtains that:
\begin{equation}\label{ReA}
F=G^2+H^2.
\end{equation}

We easily see that every representation of this type comes from a partition of the roots:

\begin{teo}\label{fgh}
 Let $F$ be a binary form of degree $2d$ and $G,H$ polynomials with real coefficients such that $F=G^2+H^2$. Then there exists a partition of its roots
$(A,\overline{A})$ such that without loss of generality $G$ is the real part and $H$ is the imaginary part of $A$.
\begin{proof}
 Let $F$ be as in the hypothesis. So $F=(G+iH)(G-iH)$ and the polynomials $G+iH$ e $G-iH$ have degree $d$. In the complex field the two polynomials are
product of $d$ linear forms, and so all their roots are also roots of $F$. If $G(\alpha)+i\cdot H(\alpha)=0$ then
$G(\overline{\alpha})-i\cdot G(\overline{\alpha})=0$. With $A=G+iH$ we conclude.
\end{proof}
\end{teo}

\begin{lemma}
Let $d \geq 2$ and $A=\prod_{i=1}^d{(x-\alpha^{(i)} y)}\in\mathbb{C}[x,y]$, with $\mathrm{Im}(\alpha^{(i)})\mathrm{Im}(\alpha^{(j)})>0$
for every $i,j$. If $A=G(x,y)+iH(x,y)$ with $G,H\in\mathbb{R}[x,y]$, then $G$ and $H$ have only real roots.
\begin{proof}
For every complex number $\gamma$ let $\gamma_1=\Ree\gamma$ and $\gamma_2=\Imm\gamma$. We proceed by induction on $d$. If $d=2$
$$
  A=(x-\alpha y)(x-\beta y)=\Big(x^2-(\alpha_1+\beta_1)xy+(\alpha_1\beta_1-\alpha_2\beta_2)y^2\Big)+i\Big(-(\alpha_2+\beta_2)xy+(\alpha_1\beta_2+\alpha_2\beta_1)y^2\Big).
$$
It is easy to see that the imaginary part of $A$ has only real roots, while the real part has only real roots if and only if its discriminant
\begin{equation}
 \Delta=\Big(\alpha_1-\beta_1\Big)^2+4\alpha_2\beta_2
\end{equation}
is nonnegative. The hypothesis about the sign of $\alpha_2\cdot\beta_2$ concludes the first part. Suppose now that $d\geq 3$. If $A=G+iH$, then
$A_x=G_x+iH_x$ and $A_y=G_y+iH_y$ have degree $d-1$; moreover, by Gauss-Lucas Theorem, the roots of $A_x$ and $A_y$ lie in the convex hull of the
set of roots of $A$, and so also the imaginary parts of the roots of $A_x,A_y$ have the same sign. The same applies for any substitution
$\tilde{A}(t)=A(at+c,bt+d)$ (because $\text{SL}(2,\mathbb{R})$ acts both on the upper half plane and on the lower half-plane),
and so, for every $\mathbb{R}^2\ni(a,b)\neq(0,0)$, the imaginary parts of the roots of $aA_x+bA_y=\frac{d}{dt}\tilde{A}(t)$ have the same sign.
By induction, since
\begin{equation}
 aA_x+bA_y=\Big(aG_x+bG_y\Big)+i\Big(aH_x+bH_y\Big),
\end{equation}
the polynomials $aG_x+bG_y$ and $aH_x+bH_y$ have only real roots
for every $\mathbb{R}^2\ni(a,b)\neq(0,0)$. So, by \cite[Th.1]{causa2011maximum} we conclude that $G$ and $H$ have only real roots.
\end{proof}
\end{lemma}

We observe that the statement of \cite{causa2011maximum} we use, is closely related to Obreschkoff theorem (see for example \cite{dedieu1992obreschkoff}).
Now we are able to deduce the main result of this paper:

\begin{teo}
 For every $d \geq 2$, $\mathcal{C}(\mathbf{P}_{2,2d})=2$.
 \begin{proof}
  Let $F\in\mathbf{P}_{2,2d}$ and let $(A,\overline{A})$ be the partition of the roots of $F$ obtained by choosing in $A$ the roots
  with positive imaginary part (or negative, that is the same). Then, if $A=G+iH$, $G$ and $H$ have only real roots and $F=G^2+H^2$
  is an extremal decomposition for $F$.
 \end{proof}
\end{teo}

Recalling that every representation of a nonnegative binary form as a sum of two squares comes from a partition of the roots, and that the imaginary part of
the polynomial $A$ of any partition is a multiple of $y$, we deduce the following fact about canonical representations of nonnegative binary forms:

\begin{coroll}\label{main}
 Let $d \geq 1$ and let $F=F(x,y)$ be a nonnegative binary form of degree $2d$ whose term on $x^d$ has coefficient 1. Then there exist two binary forms
 $\mathcal{L},\mathcal{M}$ of degree $d$ and $d-1$ respectively, with only real roots, and such that $$F=\mathcal{L}^2+y^2\mathcal{M}^2.$$ 
\end{coroll}

\bibliographystyle{plain}

\bibliography{biblio}

\end{document}